\documentclass[a4paper,12pt]{amsart}
\usepackage{amssymb}
\usepackage{ifthen}
\usepackage{graphicx}

\newtheorem{thm}{Theorem}
\newtheorem{cor}{Corollary}
\newtheorem{lem}{Lemma}

\newtheorem{rem}{Remark}

\newtheorem{conj}{Conjecture}
\theoremstyle{definition}
\newtheorem{example}[equation]{Example}
\newtheorem{prob}[equation]{Problem}

\newcounter {own}
\def\theown {\thesection       .\arabic{own}}

\newenvironment{pf}[1][]{%
 \vskip 3mm
 \noindent
 \ifthenelse{\equal{#1}{}}%
  {{\slshape Proof. }}%
  {{\slshape #1.} }%
 }%
{\qed\bigskip}

\newcounter{alphabet}
\newcounter{tmp}
\newenvironment{Thm}[1][]{\refstepcounter{alphabet}%
\bigskip%
\noindent%
{\bf Theorem \Alph{alphabet}}%
\ifthenelse{\equal{#1}{}}{}{ (#1)}%
{\bf .} \itshape}{\vskip 8pt}

\makeatletter
\newcommand{\Ref}[1]{\@ifundefined{r@#1}{}{\setcounter{tmp}{\ref{#1}}\Alph{tmp}}}
\makeatother

\newenvironment{Lem}[1][]{\refstepcounter{alphabet}%
\bigskip%
\noindent%
{\bf Lemma \Alph{alphabet}}%
{\bf .} \itshape}{\vskip 8pt}

\newcommand{\IR}{{\mathbb R}}

\newcommand{\IC}{{\mathbb C}}
\newcommand{\ID}{{\mathbb D}}

\def\be{\begin{equation}}
\def\ee{\end{equation}}

\newcommand{\bee}{\begin{enumerate}}
\newcommand{\eee}{\end{enumerate}}

\newcommand{\blem}{\begin{lem}}
\newcommand{\elem}{\end{lem}}
\newcommand{\bthm}{\begin{thm}}
\newcommand{\ethm}{\end{thm}}
\newcommand{\bcor}{\begin{cor}}
\newcommand{\ecor}{\end{cor}}
\newcommand{\beg}{\begin{example}}
\newcommand{\eeg}{\end{example}}
\newcommand{\begs}{\begin{examples}}
\newcommand{\eegs}{\end{examples}}
\newcommand{\bdefe}{\begin{defin}}
\newcommand{\edefe}{\end{defin}}
\newcommand{\bprob}{\begin{prob}}
\newcommand{\eprob}{\end{prob}}
\newcommand{\bei}{\begin{itemize}}
\newcommand{\eei}{\end{itemize}}

\newcommand{\bcon}{\begin{conj}}
\newcommand{\econ}{\end{conj}}
\newcommand{\bcons}{\begin{conjs}}
\newcommand{\econs}{\end{conjs}}
\newcommand{\bprop}{\begin{propo}}
\newcommand{\eprop}{\end{propo}}
\newcommand{\br}{\begin{rem}}
\newcommand{\er}{\end{rem}}
\newcommand{\brs}{\begin{rems}}
\newcommand{\ers}{\end{rems}}
\newcommand{\bo}{\begin{obser}}
\newcommand{\eo}{\end{obser}}
\newcommand{\bos}{\begin{obsers}}
\newcommand{\eos}{\end{obsers}}
\newcommand{\bpf}{\begin{pf}}
\newcommand{\epf}{\end{pf}}
\newcommand{\ba}{\begin{array}}
\newcommand{\ea}{\end{array}}
\newcommand{\beq}{\begin{eqnarray}}
\newcommand{\beqq}{\begin{eqnarray*}}
\newcommand{\eeq}{\end{eqnarray}}
\newcommand{\eeqq}{\end{eqnarray*}}

\def\cc{\setcounter{equation}{0}   
\setcounter{figure}{0}\setcounter{table}{0}}

\begin{document}
\bibliographystyle{amsplain}
\title[On relations between the classes $\mathcal S$ and  $\mathcal U$
]{On relations between the classes $\mathcal S$ and  $\mathcal U$}

\author[M. Obradovi\'{c}]{Milutin Obradovi\'{c}}
\address{M. Obradovi\'{c},
Department of Mathematics,
Faculty of Civil Engineering, University of Belgrade,
Bulevar Kralja Aleksandra 73, 11000
Belgrade, Serbia. }
\email{obrad@grf.bg.ac.rs}

\author[S. Ponnusamy]{Saminathan Ponnusamy
$^\dagger $
}
\address{S. Ponnusamy, Indian Statistical Institute (ISI), Chennai Centre, SETS (Society
for Electronic Transactions and Security), MGR Knowledge City, CIT
Campus, Taramani, Chennai 600 113, India.}
\email{samy@isichennai.res.in, samy@iitm.ac.in}

\author[K.-J. Wirths]{Karl-Joachim Wirths}
\address{K.-J. Wirths, Institut f\"ur Analysis und Algebra,
TU Braunschweig, Pockelsstr. 14 D 38106 Braunschweig, Germany.}
\email{kjwirths@tu-bs.de}

\keywords{Analytic and univalent functions, subordination, Koebe transforms, and Grunsky inequalities \\
$^\dagger $ Corresponding author
}
\date{\today  
;  File: }

\begin{abstract}
Let ${\mathcal A}$ denote the family of all functions $f$ analytic
in the unit disk $\ID$ and satisfying the normalization
$f(0)=0= f'(0)-1$. Let $\mathcal{S}$ denote the subclass of ${\mathcal A}$ consisting of univalent functions  in  $\ID$. We consider the
subclass $\mathcal{U} $ of  $\mathcal{S}$ that is defined by the condition that for its members $f$ the condition
$$\left |\left (\frac{z}{f(z)} \right )^{2}f'(z)-1\right | < 1 ~\mbox{ for $z\in \ID$}
 $$
holds. To theses relations belong striking similarities and on the other hand big differences.
We show that some results about $\mathcal{S}$ can be improved for $\mathcal{U}$, while others cannot.
\end{abstract}

\thanks{Dedicated to Professor David Minda}

\maketitle
\pagestyle{myheadings}
\markboth{M. Obradovi\'{c},  S. Ponnusamy, and K.-J. Wirths }{On relations between the classes $\mathcal S$ and  $\mathcal U$}
\cc

\section{Introduction and Statement of Results}

Let ${\mathcal A}$ denote the family of all functions $f$ analytic
in the unit disk $\ID := \{ z\in \IC:\, |z| < 1 \}$ and  satisfying the normalization
$f(0)=0= f'(0)-1$. Let $\mathcal{S}$ denote the subclass of ${\mathcal A}$ consisting of univalent functions in  $\ID$.
We consider relationships between $\mathcal{S}$ and its subclass $\mathcal{U} $ that is defined by the condition that for its members $f$ the condition
\be\label{OPW-eq3}
\left |\left (\frac{z}{f(z)} \right )^{2}f'(z)-1\right | < 1
 ~\mbox{ for $z\in \ID$}
 \ee
holds. It has been proved in \cite{A} that $\mathcal{U} \subset \mathcal{S}$. Typical members of the class $\mathcal{U}$
are
$$\frac{z}{1\pm z},\; \frac{z}{1\pm z^2},\; \frac{z}{(1\pm z)^2},\;
\frac{z}{1\pm z+z^2}
$$
and their rotations. The class $\mathcal{U}$ and its various generalizations have been studied recently. In particular,
the class $\mathcal{U}$ is  preserved under rotation, conjugation, dilation and omitted-value transformations
but is not preserved under the square-root transformation, for example. See \cite{OPW} and the references therein.

In the present paper we consider some problems, where the solutions are identical for $\mathcal{S}$ and $\mathcal{U}$ and
some others where there exist differences.

The first problem we address is the question for the maximum radius of the circle around the origin wherein ${\rm Re} (f(z)/z)\,>\,1/2$. The solution
for the class $\mathcal{S}$ was presented in \cite{S} and \cite{W} as follows: If  $f\in \mathcal{S}$, then
\be\label{OPW-eq1}
{\rm Re}\left(\frac{f(z)}{z}\right) > \frac{1}{2}
\ee
if $|z|<\sqrt{2}-1$. This bound is best possible.


It is worth recalling that if $f\in \mathcal{S}$ is convex or starlike of order $1/2$ or $f\in \mathcal{A}$ such that the Taylor coefficients of $f$
are real and convex decreasing, then the condition \eqref{OPW-eq1} holds in the full disk $\ID$. Secondly,  since $\mathcal{U} \subsetneq \mathcal{S}$,
\eqref{OPW-eq1} holds for the class $\mathcal{U}$, too. Indeed,
the Koebe function $ k(z)= z/(1-z)^2$  belongs to the class $\mathcal{U}$ and the equation $r^{-1}k(r)=1/2$ with $r=1-\sqrt{2}$
as well as the considerations in \cite{S} show that the result \eqref{OPW-eq1} is still the best possible for the class $\mathcal{U}$.



The situation changes significantly if one considers the similar problem asking where
\be\label{eq-new1}
{\rm Re}\left ( \sqrt{\frac{f(z)}{z}}\right )>\frac{1}{2}
\ee
is valid. In 1971,  the following result was proved by Duren and Schober \cite{DS}.

\begin{Thm}\label{OPW-Th2}
For each $f \in \mathcal{S}$, the inequality \eqref{eq-new1} holds for $|z|<R$, where $R=0.835\ldots$ is the best possible radius. 
Moreover, for each $z$ in $|z|>R$, there exists an $f_0\in\mathcal{S}$ for which \eqref{eq-new1} fails to hold.
\end{Thm}

Concerning the same question for the class $\mathcal{U}$, we may recall the following result of Obradovi\'{c} proved in \cite{O}:
\be\label{OPW-eq2}
f\in \mathcal{U} \Rightarrow \frac{z}{f(z)}\prec (1-z)^{2}, ~\mbox{ i.e. }~{\rm Re}\left ( \sqrt{\frac{f(z)}{z}}\right )>\frac{1}{2}
\ee
is valid for $z\in\ID$ (see also \cite{OPW}).  Here $\prec$ denotes the usual subordination (cf. \cite{MM-book,OPW,P}).


%
In the following we will generalize the implication \eqref{OPW-eq2} for the class $\mathcal{U}_n:= \mathcal{A}_{n}\cap \mathcal{U}$,
where $\mathcal{A}_{n},\, n\geq 1$, denotes the class of  functions $f\in \mathcal{A}$ of the form
$$f(z)=z+a_{n+1}z^{n+1}+ \cdots  .
$$

\bthm\label{th3}
If $f\in \mathcal{U}_{n}$, then
\be\label{eq14}
{\rm Re} \left(\frac{f(z)}{z}\right)^{\frac{n}{2}}>\frac{1}{2} ~\mbox{ for }~z\in \ID .
\ee
\ethm

For $n=1$, it is a simple corollary to Theorem \Ref{OPW-Th2} that this stands in contrast to the situation in the class $\mathcal{S}$.
Choose the  function $f_0\in \mathcal{S}$ and the number $z_0$ as indicated in Theorem \Ref{OPW-Th2} and let $f_0(z)\,=\,zh_0(z)$.
Then we have
\[{\rm Re}\left (\sqrt{h_0(z_0)} \right )< \frac{1}{2}.
\]
Let further
\[g_0(z) = z\sqrt[n]{\frac{f_0(z^n)}{z^n}}
\]
and choose $z_1$ such that $z_1^n = z_0$, where $z_0$ is a complex number such that $|z_0|=R=0.835\ldots$ and thus, 
$|z_1|=\sqrt[n]{|z_0|}=\sqrt[n]{0.835\ldots}$. Then $g_0(z)\in \mathcal{S}\cap \mathcal{A}_{n} $ and
\[
{\rm Re} \left(\frac{g_0(z_1)}{z_1}\right)^{\frac{n}{2}}\,=\,{\rm Re}\left (\sqrt{h_0(z_0)}\right ) <\frac{1}{2}.
\]

Another item where one can see as well similarities as differences between the two classes in question is the problem of
Koebe transforms. For $f \in \mathcal{S}$, we define the Koebe transform with respect to the point $\zeta \in \ID$ as
\[g(z):= g(\zeta,z)= \frac{f\left(\frac{\zeta +z}{1+\overline{\zeta}z}\right) - f(\zeta)}{f'(\zeta)(1-|\zeta |^2)}.
\]
Then it is well known that these Koebe transforms as functions of the variable $z$ are all members of the class $\mathcal{S}$.

For the class $\mathcal{U}$ we prove.

\bthm\label{th4}
Let $f\in \mathcal{U}$. Then the Koebe transforms of $f$ with respect to any fixed $\zeta$, i.e. the functions $z\mapsto g(z) $ as above,
belong to $ \mathcal{U}$ if and only if
\be\label{e1}
\left|\frac{(\zeta -u)^2f'(\zeta )f'(u)}{(f(u)-f(\zeta))^2}- 1\right| < 1,\quad \zeta,\, u \in \ID.
\ee
\ethm
%

Remarkably, the disk with center at the origin, wherein (\ref{e1}) is satisfied for all members of the class, is the same for
the classes $ \mathcal{S}$ and $ \mathcal{U}.$ Finally, we also prove

\bthm\label{th5}
Let $f\in \mathcal{S}$ or $f\in  \mathcal{U}$. Then the inequality \eqref{e1} is satisfied for $|\zeta|, \,|u| < \sqrt{2}-1$.
The result is best possible in both cases.
\ethm

We note that it might be worthwhile to consider those functions that satisfy the
condition of Theorem \ref{th4}.

The proofs of Theorems \ref{th3}, \ref{th4} and \ref{th5} will be presented in Section \ref{sec2}.

\section{Proofs of Theorems \ref{th3}, \ref{th4} and \ref{th5}}\label{sec2}

The following lemma due to Miller and Mocanu \cite{MM78} is needed for the proof of Theorem \ref{th3}. See \cite{MM-book} for a general formulation of this lemma
via differential subordination.

\begin{Lem}{\rm \cite{MM78}}\label{MM-lem1}
Suppose that $\psi:\,\IC ^2 \rightarrow \IC$ is continuous in a domain $D$ of $\IC^2$ such that
$(1,0)\in D$,  ${\rm Re}\,\psi (1,0)>0$ and
$${\rm Re}\,\psi  (ix,y)\leq 0  ~\mbox{ for all  $(ix,y)\in D$ and $y\leq -n(1+x^2)/2$}.
$$
where $n\geq 1$. Let $p(z)=1+p_nz^n+\cdots $ be analytic in $\ID$ and $p(z)\not\equiv 1$. If $(p(z),zp'(z))\in D$ for all $z\in \ID$ and
${\rm Re}\,\psi (p(z),zp'(z))>0$ for all $z\in \ID$,
then ${\rm Re}\,p(z)>0$ in $\ID$.
\end{Lem}

%

\subsection{Proof of Theorem \ref{th3}}
For $n=1$, the result is the content of the implication \eqref{OPW-eq2}.
For $n=2$ (i.e. when $a_{2}=0$), the appropriate result is given in the paper \cite{OP} but the same may
be obtained from the proof that follows now.

Let $f\in \mathcal{U}$. Then  \eqref{OPW-eq3} holds, or equivalently
$${\rm Re} \left (2\left(\frac{f(z)}{z}\right)^{2}\frac{1}{f'(z)}-1\right)>0~\mbox{ for $z\in \ID$}.
$$
We now introduce
\be\label{eq15}
p(z)=2\left(\frac{f(z)}{z}\right)^{\frac{n}{2}}-1.
\ee
Clearly,  $p$ is analytic in $\ID$ and has the form $p(z)=1+p_{n}z^{n}+ \cdots .$ We shall apply
Lemma \Ref{MM-lem1} and prove that ${\rm Re}\, p(z)>0$ for $z\in\ID$.  From \eqref{eq15} we have
$$\frac{f(z)}{z}=\left(\frac{p(z)+1}{2}\right)^{\frac{2}{n}}
$$
and a computation gives that
\beqq
2\left(\frac{f(z)}{z}\right)^{2}\frac{1}{f'(z)}-1 &=:& \psi (p(z),z p'(z)),
\eeqq
where
\be\label{eq16}
\psi (r,s)=\frac{2n\left(\frac{r+1}{2}\right)^{\frac{2}{n}}(r+1)}{n(r+1)+2s}-1 .
\ee
According to Lemma \Ref{MM-lem1}, to prove ${\rm Re}\, p(z)>0$ in $\ID$, it suffices to show that 
\be\label{MM-eq2}
{\rm Re\,} \psi (ix,y)\leq 0 ~\mbox{ for all reals $x,y$  with $y\leq -n(1+x^2)/2$.}
\ee
It follows that
\be\label{eq17}
{\rm Re\,} \psi (ix,y)={\rm Re}\,\frac{2n\left(\frac{ix+1}{2}\right)^{\frac{2}{n}}(ix+1)}{n(ix+1)+2y}-1 .
\ee
We may use the representation $ix+1= (\sqrt{1+x^2})\, e^{i\varphi}$, $|\varphi|< \frac{\pi}{2}$,
where
\be\label{eq19}
\cos\varphi= \frac{1}{\sqrt{1+x^2}},\,\,\sin\varphi=\frac{x}{\sqrt{1+x^2}} ~\mbox{ and } ~\tan\varphi =x.
\ee
Clearly, $x \sin\varphi\geq 0 $ and, since $n\geq 2$, $x \sin\left(\frac{2}{n}\varphi\right)\geq 0 $.

By using \eqref{eq17} and \eqref{eq19}, after some simple transformations, we obtain that
$${\rm Re\,} \psi (ix,y)=\frac{S-T}{(n+2y)^{2}+n^{2}x^{2}},
$$
where
$$S=2n\left(\frac{1+x^{2}}{4}\right)^{\frac{1}{n}}\cos\left(\frac{2}{n}\varphi\right)(n+2y +nx^{2})
$$
and
\beq
T 
&=& 4y^{2}+ 4n\left(1+\left(\frac{1+x^{2}}{4}\right)^{\frac{1}{n}}x\sin\left(\frac{2}{n}\varphi\right)\right)y +n^{2}(1+x^{2}).
\label{eq20}
\eeq
Clearly  $S\leq 0$ for all $n\geq 2$ and for  $y\leq -(n/2)(1+x^{2})$.
Thus, we also need to prove that $T:=T(y)\geq 0$ for $n\geq 2$ and for all $x\in\IR$ and  $y\leq -(n/2)(1+x^{2})$.
The function $T(y)$  has its minimum value at the point
$$ y_{0}=-\frac{n}{2}\left(1+\left(\frac{1+x^{2}}{4}\right)^{\frac{1}{n}}x\sin\left(\frac{2}{n}\varphi\right)\right)
$$
so that $T(y)\geq T(y_0)$. Since
$$
\left(\frac{1+x^{2}}{4}\right)^{\frac{1}{n}}x\sin\left(\frac{2}{n}\varphi\right)
\leq \left(\frac{1+x^{2}}{4}\right)^{\frac{1}{n}}x\sin \varphi
=\frac{x^{2}}{2}\left(\frac{4}{1+x^{2}}\right)^{\frac{1}{2}-\frac{1}{n}}
\leq  x^{2},
$$
for $n\geq 2$,  we easily conclude that $- (n/2)(1+x^{2})\leq y_{0}$. As $T(y)$ is decreasing when $y\leq y_0$, it is enough to prove that
\be\label{eq21}
T\left(-\frac{n}{2} (1+x^{2})\right)= n^{2}(1+x^{2})\left [x^{2}
-2\left(\frac{1+x^{2}}{4}\right)^{\frac{1}{n}}x\sin\left(\frac{2}{n}\varphi\right)\right ]\geq 0
\ee
for all  $x\in \mathbb{R}$ and $n\geq 2$.
Since, by the previous consideration, $x\sin\left(\frac{2}{n}\varphi\right)\geq0$, we can suppose that
$x \geq 0$ and $0\leq \varphi<\frac{\pi}{2}$.
In view of this observation, proving the inequality \eqref{eq21} is equivalent to proving the inequality
\be\label{eq22}
\sin\left(\frac{2}{n}\varphi\right)\leq\frac{x}{2} \left(\frac{4}{1+x^{2}}\right)^{\frac{1}{n}} ~\mbox{ for $x\geq 0$, $0\leq \varphi<\pi/2$, and $n\geq 2$}.
\ee
For $n=2$, we have equality in \eqref{eq22} (by using \eqref{eq19}).
Again, from \eqref{eq19}, we obtain  that $\sin^{2}\varphi=x^{2}/(1+x^{2})$ and $x=\tan\varphi$, and
thus the inequality \eqref{eq22} is equivalent to the inequality
\be\label{eq23}
g(\varphi)\geq g(0)=0 ~\mbox{ for $0\leq \varphi<\pi/2$ and $n\geq 2$},
\ee
where
$$g(\varphi)=(2\cos\varphi)^{\frac{2}{n}-1}\sin\varphi-\sin\left(\frac{2}{n}\varphi\right).
$$
We find that
$$g'(\varphi)=2(2\cos\varphi)^{\frac{2}{n}-2}\left(1-\frac{2}{n}\sin^{2}\varphi\right)-\frac{2}{n}\cos\left(\frac{2}{n}\varphi\right)
$$
and thus,
$$g'(0)=2\left(\frac{1}{2^{2-\frac{2}{n}}}-\frac{1}{n}\right)>0~\mbox{ for $n\geq 3$.}
$$
Also, a computation gives that
$$g''(\varphi)=8(2\cos\varphi)^{\frac{2}{n}-3}\left(\left(1-\frac{3}{n}\right)\sin \varphi
+\frac{2}{n^{2}}\sin^{3}\varphi\right)+\frac{4}{n^{2}}\sin \left(\frac{2}{n}\varphi\right)\geq 0~\mbox{ for $n\geq 3$.}
$$
It means that the function $g'$ is an increasing function of $\varphi$ and this gives
$$g'(\varphi)\geq g'(0)>0  ~\mbox{for }~ 0\leq\varphi<\pi/2
$$
which in turn implies that the function $g(\varphi)$ is also increasing for $0\leq\varphi<\pi/2$ and  hence, \eqref{eq23} holds.
This means that \eqref{MM-eq2} holds and hence, by Lemma \Ref{MM-lem1}, it follows that
${\rm Re}\,p(z)>0$ in $\ID$.  The proof of the theorem is complete.
\qed

\vspace{8pt}

For the proof of Theorem \ref{th4}, we need the following lemma, which might have been known in the literature.
Since we were not able to find an apt reference we give the proof for this theorem.
We want to emphasize here that the functions considered in this lemma are neither conformal maps nor harmonic functions.

\blem\label{th6} Let for $z\in \ID$,
\[ u(\zeta)=\frac{z+\zeta}{1+\overline{\zeta}z}.
\]
Then $u:\ID\to\ID$ and $u:\overline{\ID} \to \overline{\ID}$ are bijective.
\elem
\bpf
The injectivity is easily derived from
\[ u(\zeta_1)-u(\zeta_2)
\,=\,\frac{(\zeta_1-\zeta_2)(1+\overline{\zeta_1}z) +(\overline{\zeta_1}-\overline{\zeta_2})(\zeta_1z+z^2)}{(1+\overline{\zeta_1}z)(1+\overline{\zeta_1}z)}\]
If this difference equals zero and $\zeta_1\,\neq \zeta_2$, then
\[ |1+\overline{\zeta_1}z|=|\zeta_1z+z^2|, \quad (1-|z|^2)(1+|z|^2 +2{\rm Re}\, (\overline{\zeta_1}z)=0,
\]
which is impossible for $z\in \ID$.

Further the functional determinant
\[ \left|\frac{\partial u}{\partial \zeta}\right|^2\,-\left|\frac{\partial u}{\partial \overline{\zeta}}\right|^2
= \frac{|1+\overline{\zeta}z|^2 - |\zeta z+z^2|^2}{|1+\overline{\zeta}z|^4}
\]
does not equal zero for $\zeta \in \ID$. Hence $u(\ID)$ is open. Further it is easily seen that $u(\ID)\subset \ID$ and $u(\partial \ID)=\partial \ID.$

Now assume that $\ID\setminus u(\ID)$ is non-void and not open. Then there exists a point $p\in \ID\setminus u(\ID)$ and a
sequence $\{\zeta_n\}_{n\geq 1}$ in $\ID,$ such that
\[ p=\lim_{n\to\infty}u(\zeta_n).
\]
Let $\{\zeta_{n_k}\}_{k\geq 1}$ be a convergent subsequence of $\{\zeta_n\}$. Because of the continuity of $u$ and $u(\partial \ID)=\partial \ID$,
we have $\lim_{k\to\infty}\zeta_{n_k}\,=\,w \in \ID$ and
\[ p=\lim_{k\to\infty}u(\zeta_{n_k})=u(w).
\]
Hence, $\ID\setminus u(\ID)$ is void or open. The second possibility contradicts the connectivity of $\ID$.
Together with the above this proves the assertions.
\epf

\subsection{Proof of Theorem \ref{th4}} Let $f$ belong to $ \mathcal{U}$ and for fixed $\zeta\in\ID$, consider its Koebe transforms $g(z)$
with respect to $\zeta$ given by
$$g(z):= g(\zeta,z)= \frac{f\left( u(\zeta)\right) - f(\zeta)}{f'(\zeta)(1-|\zeta |^2)}, \quad  u(\zeta)=\frac{z+\zeta}{1+\overline{\zeta}z}.
$$
If all Koebe transforms of $f$ belong to $ \mathcal{U}$, then by \eqref{OPW-eq3} we have
\[\left|\frac{z^2f'(\zeta)(1-|\zeta|^2)^2}{(f(u(\zeta))-f(\zeta ))^2}\frac{f'(u(\zeta))}{(1+\overline{\zeta}z)^2} - 1\right|=
\left|\frac{(\zeta-u(\zeta))^2f'(\zeta)f'(u(\zeta))}{(f(u(\zeta))-f(\zeta ))^2} - 1\right|<1
\]
for all $u,\zeta \in \ID$. According to Lemma \ref{th6} this proves the necessity of the above condition.
The sufficiency can be proved similarly. \hfill $\Box$

\subsection{Proof of Theorem \ref{th5}}
Let $f\in\mathcal{S}$ and let
\be\label{eq38}
\log \frac{f(z)-f(u)}{z-u}=\sum_{n,m=0}^{\infty}d_{n,m}z^{n}u^{m}
\ee
The coefficients  $d_{n,m}$ are called the Grunsky coefficients of the function $f$.  From \eqref{eq38}, after differentiations
with respect to $z$ and $u$, we have
$$\frac{f'(z)f'(u)}{(f(z)-f(u))^{2}}-\frac{1}{(z-u)^{2}}=\sum_{n,m=1}^{\infty}nmd_{n,m}z^{n-1}u^{m-1}$$
 and from here
\be\label{eq39}
\frac{(z-u)^{2}f'(z)f'(u)}{(f(z)-f(u))^{2}}-1=(z-u)^{2}\sum_{n,m=1}^{\infty}nmd_{n,m}z^{n-1}u^{m-1}.
\ee

By using  Grunsky's inequalities (see \cite[p.~62]{P})
\be\label{eq 40}
\sum_{n=1}^{\infty}n\left|\sum_{m=1}^{\infty}d_{nm}x_{m}\right|^{2}\leq \sum_{n=1}^{\infty}\frac{|x_{n}|^{2}}{n},
\ee
if the last series converges and for arbitrary $x_{n},\,\,n=1,2,\ldots $
(We note that Grunsky's inequality usually is stated with the functions from the class $\Sigma $, but it is easy
to prove that Grunsky's coefficients for the functions $\log \frac{f(z)-f(u)}{z-u}$ and
$\log \frac{F(z^{-1})-F(u^{-1})}{z^{-1}-u^{-1}}$, where $F(\zeta)=\frac{1}{f(1/\zeta)}\in \Sigma $ for
$f\in \mathcal{S}$, are  the same for $n,m\geq 1.$) we can obtain that
\beqq
\left|\sum_{n,m=1}^{\infty}nmd_{n,m}z^{n-1}u^{m-1}\right|
&=& \left|\sum_{n=1}^{\infty}\sqrt{n}z^{n-1}\sqrt{n}\sum_{m=1}^{\infty}d_{n,m}mu^{m-1}\right|\\
&\leq& \left(\sum_{n=1}^{\infty}n|z|^{2(n-1)}\right)^{\frac{1}{2}}
\left(\sum_{n=1}^{\infty}n\left|\sum_{m=1}^{\infty}d_{n,m}mu^{m-1}\right|^{2}\right)^{\frac{1}{2}}\\
&\leq& \frac{1}{1-|z|^{2}}\left(\sum_{n=1}^{\infty}n|u|^{2(n-1)}\right)^{\frac{1}{2}}\\
&=&\frac{1}{(1-|z|^{2})(1-|u|^{2})}.
\eeqq
From this and \eqref{eq39} we finally have
$$\left| \frac{(z-u)^{2}f'(z)f'(u)}{(f(z)-f(u))^{2}}-1\right|\leq \frac{|z-u|^{2}}{(1-|z|^{2})(1-|u|^{2})}
\leq \left(\frac{2r}{1-r^{2}}\right)^{2}<1,
$$
since $|z|,\, |u|\leq r<\sqrt{2}-1.$

To prove that this result is sharp for $\mathcal{U}$ and $\mathcal{S}$ we consider the Koebe function $k$ that belongs to both classes.
A simple calculation reveals that for $f=k$, (\ref{e1}) becomes
$$\left | \frac{u- \zeta }{1-u\zeta }\right |<1.
$$
For $\zeta =(\sqrt{2}-1)i$ and  $u =-(\sqrt{2}-1)i$,
$$\left | \frac{u- \zeta }{1-u\zeta }\right |=1.
$$
This implies that $\sqrt{2}-1$ is best possible.
\hfill $\Box$

\section{Concluding remarks}
A natural question is the following: Are all functions
$$f(z)=z\prod _{k=1}^n(1-e^{i\theta _k}z)^{-\alpha _k},
~\mbox{ $\theta _k\in \IR$, $\alpha  _{k}\geq 0$ and $\sum\limits_{k=1}^{n}\alpha  _{k}=2$},
$$
in the class $\mathcal{U}$? The answer is no as the function $f_0(z)=z(1-z^3)^{-2/3}$ demonstrates. Note that
$$\frac{z}{f_0(z)}= (1-z)^{2/3} (1+e^{-i\pi /3}z)^{2/3} (1+e^{i\pi /3}z)^{2/3}.
$$

Moreover, if $f \in\mathcal{S}$ then $r^{-1}f(rz)\in \mathcal{U}$ for $0<r\leq1/\sqrt{2}$ and the result is sharp. See \cite{obpo-2005}.
Furthermore, the family $\mathcal{U}$ is not a subset of the class $\mathcal{S}^{\star}$ of univalent starlike functions in the unit disk $\ID$.
In fact,  consider the function
$$f_1(z)=\frac{z}{1+\frac{1}{2}z+\frac{1}{2}z^3}.
$$
Then it is easy to see that $f_1\in {\mathcal U}$. On the other hand,
$$\frac{zf_1'(z)}{f_1(z)}=\frac{1-z^3}{1+\frac{1}{2}z+\frac{1}{2}z^3}
$$
and at the boundary point $z_0= (-1+i)/\sqrt 2$, we have
$$\frac{z_0f_1'(z_0)}{f_1(z_0)} =\frac{2-2\sqrt{2}}{3} +\frac{1-2\sqrt 2}{3}i
$$
which gives that ${\rm Re}\, \{z_0f_1'(z_0)/f_1(z_0)\}<0$.
Consequently, there are points in the unit disk $|z|<1$ for which
${\rm Re}\, \{zf_1'(z)/f_1(z)\}<0$ which shows that the function
$f_{1}$ is not starlike in $\ID$.

\subsection*{Acknowledgements}
The authors thank the referee for his/her careful reading and many useful comments.
The work of the first author was supported by MNZZS Grant, No. ON174017, Serbia.
The second author is on leave from the IIT Madras.

\end{document}